\documentclass[11pt]{article}
\usepackage{amsmath,amsthm,amssymb}
\usepackage[dvips]{graphicx}

\newcommand{\C}{\mathbb C}
\newcommand{\ii}{{\operatorname{i}}}
\newcommand{\sn}{\operatorname{sn}}
\newcommand{\cn}{\operatorname{cn}}
\newcommand{\dn}{\operatorname{dn}}

\begin{document}


\title{Inequalities for the Jacobian Elliptic Functions with Complex Modulus\footnote{published in: Journal of Mathematical Inequalities {\bf 6} (2012), 91--94.}}
\author{Klaus Schiefermayr\footnote{University of Applied Sciences Upper Austria, School of Engineering and Environmental Sciences, Stelzhamerstrasse\,23, 4600 Wels, Austria, \textsc{klaus.schiefermayr@fh-wels.at}}}
\date{}
\maketitle

\theoremstyle{plain}
\newtheorem{theorem}{Theorem}
\theoremstyle{definition}
\newtheorem*{remark}{Remark}

\begin{abstract}
Despite the fact that there is a huge amount on papers and books devoted to the theory of Jacobian elliptic functions, very little is known when the modulus $k$ of these functions lies outside the unit interval $[0,1]$. In this note, we prove some simple inequalities for the absolute value of Jacobian elliptic functions with complex modulus.
\end{abstract}



\section{Introduction and Main Result}


Consider the Jacobian elliptic functions $\sn(z,k)$, $\cn(z,k)$, $\dn(z,k)$ with complex parameter (modulus) $k\in\C$. Starting with the works of Jacobi in the 1820's until now, there exists a huge amount on papers and books devoted to the theory of Jacobian elliptic functions, see, e.g., \cite{NIST}, \cite{Lawden}, or \cite{WW}. Almost all of these contributions are restricted to the case when the modulus $k$ is in the unit interval $[0,1]$. Exceptions are, e.g., the articles of Walker \cite{Walker-2003}, \cite{Walker-2009}. Since the functions $\sn(z,k)$, $\cn(z,k)$, and $\dn(z,k)$ depend on $k^2$ rather than $k$, we shall use $m=k^2$ as parameter but use the same notation $\sn(z,m)$, $\cn(z,m)$, and $\dn(z,m)$.


There seems to exist no estimates for the absolute value of $\sn(z,m)$, $\cn(z,m)$, and $\dn(z,m)$ in terms of elementary functions. In this paper, we give such estimates.


\begin{theorem}\label{Thm}
Let $m\in\C$, $|m|\leq1$, and $z\in\C$, $|z|\leq{R}<\frac{\pi}{2}$. Then the inequalities
\begin{align}
|\sn(z,m)|&\leq\frac{\sn(|z|,m_1)}{\cn(|z|,m_1)}\leq\tan|z|,\label{sn-ineq}\\
|\cn(z,m)|&\leq\frac{1}{\cn(|z|,m_1)}\leq\frac{1}{\cos|z|},\label{cn-ineq}\\
|\dn(z,m)|&\leq\frac{\dn(|z|,m_1)}{\cn(|z|,m_1)}\leq\frac{1}{\cos|z|},\label{dn-ineq}
\end{align}
hold, where $m_1:=1-|m|\in[0,1]$. For all inequalities we have equality if $z=0$ or if $z=\ii{y}$, $-R\leq{y}\leq{R}$, and $m=1$. Further, in the first inequalities of \eqref{sn-ineq}--\eqref{dn-ineq}, equality is attained if $z=\ii{y}$, $-R\leq{y}\leq{R}$, and $m=|m|$.
\end{theorem}

\section{Proof}


\begin{proof}[\bf Proof of Theorem\,\ref{Thm}]
Starting point for the proof are the Taylor expansions of $\sn(z,m)$, $\cn(z,m)$, and $\dn(z,m)$, respectively, given by \cite[Eq.\,(3.6)]{Walker-2003}
\begin{align}
\sn(z,m)&=\sum_{n=0}^{\infty}(-1)^{n}s_n(m)\,\frac{z^{2n+1}}{(2n+1)!},\label{sn-sum}\\
\cn(z,m)&=\sum_{n=0}^{\infty}(-1)^{n}c_n(m)\,\frac{z^{2n}}{(2n)!},\label{cn-sum}\\
\dn(z,m)&=\sum_{n=0}^{\infty}(-1)^{n}d_n(m)\,\frac{z^{2n}}{(2n)!},\label{dn-sum}
\end{align}
where $s_n(m)$, $c_n(m)$, and $d_n(m)$ are polynomials with positive integer coefficients. Thus, we have the inequalities
\begin{align}
|s_n(m)|&\leq{s}_n(|m|)\leq{s}_n(1),\label{sn}\\
|c_n(m)|&\leq{c}_n(|m|)\leq{c}_n(1),\label{cn}\\
|d_n(m)|&\leq{d}_n(|m|)\leq{d}_n(1).\label{dn}
\end{align}
Note that each of the power series \eqref{sn}--\eqref{dn} are absolutely convergent for $|m|\leq1$ and $|z|<\frac{\pi}{2}$, see \cite[Thm.\,3.2]{Walker-2003}. Especially,
\begin{align}
\sn(\ii{z},m)&=\ii\sum_{n=0}^{\infty}s_n(m)\,\frac{z^{2n+1}}{(2n+1)!},\label{sni}\\
\cn(\ii{z},m)&=\sum_{n=0}^{\infty}c_n(m)\,\frac{z^{2n}}{(2n)!},\label{cni}\\
\dn(\ii{z},m)&=\sum_{n=0}^{\infty}d_n(m)\,\frac{z^{2n}}{(2n)!}.\label{dni}
\end{align}
If we put $m=1$ in formulae \eqref{sni}--\eqref{dni}, we get (note that $c_n(1)=d_n(1)$)
\begin{align}
\tanh(\ii{z})&=\ii\sum_{n=0}^{\infty}s_n(1)\,\frac{z^{2n+1}}{(2n+1)!},\label{tanh}\\
\frac{1}{\cosh(\ii{z})}&=\sum_{n=0}^{\infty}c_n(1)\,\frac{z^{2n}}{(2n)!}=\sum_{n=0}^{\infty}d_n(1)\,\frac{z^{2n}}{(2n)!}.\label{cosh}
\end{align}
Hence, the inequalities
\begin{align*}
|\sn(z,m)|&\leq\sum_{n=0}^{\infty}|s_n(m)|\,\frac{|z|^{2n+1}}{(2n+1)!}\qquad\text{by}~\eqref{sn-sum}\\
&\leq\sum_{n=0}^{\infty}s_n(|m|)\,\frac{|z|^{2n+1}}{(2n+1)!}=:(*_1)\qquad\text{by}~\eqref{sn}\\
&=\frac{1}{\ii}\,\sn(\ii|z|,|m|)\qquad\text{by}~\eqref{sni}\\
&=\frac{\sn(|z|,1-|m|)}{\cn(|z|,1-|m|)}\qquad\text{by~\cite[Eq.\,(2.6.12)]{Lawden}}
\end{align*}
and
\begin{align*}
|\cn(z,m)|&\leq\sum_{n=0}^{\infty}|c_n(m)|\,\frac{|z|^{2n}}{(2n)!}\qquad\text{by}~\eqref{cn-sum}\\
&\leq\sum_{n=0}^{\infty}c_n(|m|)\,\frac{|z|^{2n}}{(2n)!}=:(*_2)\qquad\text{by}~\eqref{cn}\\
&=\cn(\ii|z|,|m|)\qquad\text{by}~\eqref{cni}\\
&=\frac{1}{\cn(|z|,1-|m|)}\qquad\text{by~\cite[Eq.\,(2.6.12)]{Lawden}}
\end{align*}
and
\begin{align*}
|\dn(z,m)|&\leq\sum_{n=0}^{\infty}|d_n(m)|\,\frac{|z|^{2n}}{(2n)!}\qquad\text{by}~\eqref{dn-sum}\\
&\leq\sum_{n=0}^{\infty}d_n(|m|)\,\frac{|z|^{2n}}{(2n)!}=:(*_3)\qquad \text{by \eqref{dn}}\\
&=\dn(\ii|z|,|m|)\qquad\text{by}~\eqref{dni}\\
&=\frac{\dn(|z|,1-|m|)}{\cn(|z|,1-|m|)}\qquad\text{by~\cite[Eq.\,(2.6.12)]{Lawden}}
\end{align*}
hold and we have proved the first inequalities of \eqref{sn-ineq}--\eqref{dn-ineq}, respectively. Analogously,
\begin{align*}
(*_1)&\leq\sum_{n=0}^{\infty}s_n(1)\,\frac{|z|^{2n+1}}{(2n+1)!}\qquad\text{by}~\eqref{cn}\\
&=\frac{1}{\ii}\,\tanh(\ii|z|)\qquad\text{by}~\eqref{tanh}\\
&=\tan|z|
\end{align*}
and
\begin{align*}
(*_2),(*_3)&\leq\sum_{n=0}^{\infty}c_n(1)\,\frac{|z|^{2n}}{(2n)!}\qquad\text{by}~\eqref{cn}\\
&=\frac{1}{\cosh(\ii|z|)}\qquad\text{by}~\eqref{cosh}\\
&=\frac{1}{\cos|z|}
\end{align*}
holds.
\end{proof}


\begin{remark}
Another (much more complicated) possibility for proving the second inequalities of \eqref{sn-ineq}--\eqref{dn-ineq}, respectively, is the following: For fixed $u\in[0,\frac{\pi}{2}]$, consider the functions
\[
f_1(m_1):=\frac{\sn(u,m_1)}{\cn(u,m_1)},\,f_2(m_1):=\frac{1}{\cn(u,m_1)},\,f_3(m_1):=\frac{\dn(u,m_1)}{\cn(u,m_1)}
\]
for $0\leq{m}_1\leq1$. Then one has to prove that these functions are strictly monotone decreasing in $m_1$, i.e.\ $f_j'(m_1)<0$, $j=1,2,3$, using formulae (710.54), (710.57), and (710.60) of \cite{BF}.
\end{remark}


\bibliographystyle{amsplain}
\bibliography{InequalitySn}

\providecommand{\bysame}{\leavevmode\hbox to3em{\hrulefill}\thinspace}
\providecommand{\MR}{\relax\ifhmode\unskip\space\fi MR }
\providecommand{\MRhref}[2]{%
  \href{http://www.ams.org/mathscinet-getitem?mr=#1}{#2}
}
\providecommand{\href}[2]{#2}
\begin{thebibliography}{1}

\bibitem{BF}
P.F. Byrd and M.D. Friedman, \emph{Handbook of elliptic integrals for engineers
  and scientists}, Springer, 1971.

\bibitem{Lawden}
D.F. Lawden, \emph{Elliptic functions and applications}, Springer, 1989.

\bibitem{NIST}
F.W.J. Olver, D.W. Lozier, R.F. Boisvert, and C.W. Clark (eds.), \emph{N{IST}
  handbook of mathematical functions}, U.S. Department of Commerce National
  Institute of Standards and Technology, Washington, DC, 2010.

\bibitem{Walker-2003}
P.~Walker, \emph{The analyticity of {J}acobian functions with respect to the
  parameter {$k$}}, R. Soc. Lond. Proc. Ser. A Math. Phys. Eng. Sci.
  \textbf{459} (2003), 2569--2574.

\bibitem{Walker-2009}
P.L. Walker, \emph{The distribution of the zeros of {J}acobian elliptic
  functions with respect to the parameter {$k$}}, Comput. Methods Funct. Theory
  \textbf{9} (2009), 579--591.

\bibitem{WW}
E.T. Whittaker and G.N. Watson, \emph{A course of modern analysis}, Cambridge
  University Press, 1962.

\end{thebibliography}

\end{document}